%% file: main_arXiv.tex
\documentclass{article}
\usepackage[utf8]{inputenc}
\usepackage{url}
\usepackage[hidelinks]{hyperref}
\usepackage{amsmath}
\usepackage{amssymb}
\usepackage{mathrsfs}
\usepackage{mathtools}
\usepackage{stmaryrd}   
\usepackage[round, sort,comma,authoryear]{natbib} 
\usepackage{graphicx}
\usepackage{float}
\usepackage[framemethod=TikZ]{mdframed}
\usepackage{footmisc}
\usepackage{algorithm} 
\usepackage{algpseudocode}

\usepackage[format=plain,
            labelfont={bf,it},
            textfont=it]{caption}
\usepackage{subcaption}
\usepackage{makecell}

\usepackage{enumitem}
\newcounter{descriptcount}

\usepackage{multirow,tabularx}
\usepackage[title]{appendix}
\usepackage{csquotes}

\usepackage{orcidlink}

\usepackage{PRIMEarxiv}
\def\keywordname{{\bfseries \emph{Keywords}}}%
\def\keywords#1{\par\addvspace\medskipamount{\rightskip=0pt plus1cm
\def\and{\ifhmode\unskip\nobreak\fi\ $\cdot$
}\noindent\keywordname\enspace\ignorespaces#1\par}}
\SetBlockThreshold{2}

\DeclareFontFamily{OT1}{pzc}{}
\DeclareFontShape{OT1}{pzc}{m}{it}{<-> s * [1.10] pzcmi7t}{}
\DeclareMathAlphabet{\mathpzc}{OT1}{pzc}{m}{it}

\title{Scope Restriction for Scalable Real-Time Railway Rescheduling: An Exploratory Study}

\author{
Erik Nygren\thanks{Work done while at Swiss Federal Railways (SBB).}\\
\url{erik@nygren.ch},  
\And
Christian Eichenberger\footnotemark[1]\hspace{3pt}{    }\orcidlink{0000-0003-4460-5194}\\
\url{christianmarkus@oakmountain.ch}, 
\And
Emma Frejinger \orcidlink{0000-0003-1930-607X}\\
\url{emma.frejinger@unmontreal.ca}
}

\date{\today}

\DeclarePairedDelimiterX{\Set}[1]\{\}{%
  \, #1 \,
}

\begin{document}

\maketitle

\input{main}

\end{document}

%% file: main.tex
\begin{abstract}
With the aim to stimulate future research, we describe an exploratory study of a railway rescheduling problem. A widely used approach in practice and state of the art is to decompose these complex problems by geographical scope. Instead, we propose defining a \emph{core problem} that restricts a rescheduling problem in response to a disturbance to only trains that need to be rescheduled, hence restricting the scope in both time and space. In this context, the difficulty resides in defining a \emph{scoper} that can predict a subset of train services that will be affected by a given disturbance. We report preliminary results using the Flatland simulation environment that highlights the potential and challenges of this idea. 
We provide an extensible playground open-source implementation based on the Flatland railway environment and Answer-Set Programming.

\end{abstract} 

\vspace{0.5cm}
\\


%


\section{Introduction}\label{sec:introduction}

Railways face a multitude of planning problems, the most important ones can be categorized as routing or scheduling problems \citep{CordEtAl98}. Train scheduling, also known as train timetabling, is a tactical planning problem that has been extensively studied in the literature \citep[e.g.,][]{CaprEtAl02,CaccEtAl15}. A schedule is in place over a certain period of time (e.g., a season) and often only minor changes are made over different periods. These schedules should be robust to minor perturbations. 
More important perturbations referred to as disturbances \citep[e.g.,][]{LiuDess19, CaccToth12} may induce primary delay on directly affected trains which, in turn, can propagate to other trains in the network (secondary delay). Disturbances can be resolved by, e.g., changing the impacted trains' schedule or routes. 
Severe perturbations -- referred to as disruptions -- caused by, for example, extreme weather, infrastructure or equipment failures may in addition require rescheduling of rolling stock or crew. 

The focus of our work lies on real-time railway traffic management (aka traffic control, train dispatching or train rescheduling). The purpose is to identify mitigating strategies in short computing time so as to minimize the propagation caused by disturbances and recover the original schedule -- henceforth real-time rescheduling (or simply rescheduling).

There is an extensive literature on this topic and an exhaustive literature is out of the scope of this paper. Instead, we refer the reader to \cite{Luan19} for a recent review as well as surveys \citep[e.g.,][]{Lusby2011, AhujaEtAl05}. Real-time railway traffic management is challenging because (i) the computing time budget is very restricted, and (ii) controlling trains requires a detailed representation of the system (so-called microscopic level) that considers block sections and signals. For these reasons, it is not possible to solve the full traffic management problem for real railway networks in real time at a microscopic level.
Hence the problem is typically decomposed. 

\cite{LuanEtAl18} describe three decomposition methods of a time-space graph used in the literature: geography-based, train-based and time interval based. However, decomposition introduces a challenging coordination problem. Consider, for example, a decomposition according to geographical scope into different areas. In this case the decisions between the different areas need to be coordinated, e.g., through conditions at the boarders. \cite{Corman2012} study this coordination problem. A difficult challenge resides in striking the right balance between the size of the areas and the difficulty of the coordination problem. In particular the areas of interest can change depending on the nature and location of the disturbance. \cite{VanThiEtAl18} focus on conflict prevention and they propose a heuristic to dynamically define so-called impact zones.

In this work, we outline ideas for a solution approach that consists of a problem scope reduction step before solving a restricted traffic rescheduling problem. We conjecture that any rescheduling problem arising within the entire system can be reduced to a much smaller problem that specifies which subset of trains and which parts of the scheduled trainruns might be re-scheduled (adjusted). We refer to this as the \emph{core problem}. The definition of the core problem depends on the perturbation and can vary in both geographical and time dimensions. The objective is to dynamically reduce the scope of the full rescheduling problem so that it can be solved within the limited computing time budget while reducing the need for coordination (or avoiding it all together). This idea is closely related dynamic impact zones in~\cite{VanThiEtAl18}. However, the core problem is a more general concept and does not necessarily correspond to a given zone.

The following observations on the rescheduling problem constitute the rationale for our work.  First, the infrastructure (railway network) typically does not change over time. Second, train schedules are slow changing (often only minor changes are made from one schedule to the next) and a given schedule is in place for an extended period of time. Third, in practice, the problem is solved by experienced humans. 
The rescheduling problem hence arises in highly structured situations where experience from the past is valuable to find good-quality solutions. 
Therefore, we hypothesize that machine learning (ML) can be used to efficiently restrict the rescheduling problem scope.

The methodology can be viewed as a pre-processing phase that is somewhat generic to the optimization model and solution approach of the rescheduling problem. For a given disturbance, the goal is to define the core problem, i.e., the sub-problem affected by the perturbation. In our context the scope consists of trains' passing times and routes. The rescheduling problem is then solved keeping the solution fixed for all variables outside the scope of the core problem. This hence results in a heuristic solution approach. Assuming access to historical data or a simulator, we propose to learn from data how to restrict the scope to the core problem for any given perturbation. Based on the aforementioned observations, the scope restriction problem is a prime candidate for ML. However, in this first exploratory work whose purpose is to assess the potential of the approach, we resort to simple heuristics. 

The motivation for our work is supported by the findings in \cite{FiscMona17} in the domain of operations research (OR). They investigate if general-purpose Mixed-Integer Linear Programming solvers can successfully solve real-time rescheduling problems in an effective way, as opposed to using tailored heuristics. They show that thanks to simple pre-processing heuristics they can solve most of the benchmark instances within the given time budget. There are only two pre-processing operations: heuristic bound tightening and heuristic variable fixing. The idea we put forward here can be viewed as heuristic variable fixing. 

Further support comes from \cite{li_learning_2021} in the domain of ML, who address scalability of combinatorial optimization problems by learning to identify smaller sub-problems which can be readily solved by existing methods in a black-box manner.
In contrast to our setting, they start from a feasible solution and put emphasis on an iterative solution improvement, whereas we start from an ``almost feasible'' solution corrupted only by a disturbance and do not consider the iterative setting.
In addition, they also present a Transformer-based architecture for learning sub-problem selection, whereas we only show that speed-up is possible in the context of railway rescheduling problems.
In the terminology of \cite{bengio_machine_2020}, our approach  puts forward \emph{learning to configure algorithms}, i.e. augmenting an operation research algorithm with valuable pieces of information (i.e. the scope restriction).

The contributions of this paper are:
\begin{itemize}
    \item We introduce the scope restriction problem in railway re-scheduling, whose goal is to identify a core sub-problem in space and time, which can be be passed to an existing (general-purpose) solver.
    \item We assess the potential of scope restriction in an exploratory numerical study using a simulator  \citep[Flatland toolbox,][]{mohanty2020flatlandrl} and one specific optimization model and solver. Our exploratory results show that a significant computation time reduction for railway re-scheduling can be achieved when the problem scope can be reduced reliably.
    \item The main objective with sharing this preliminary work is to encourage future research on the topic. We contribute to such work by providing and describing an extensible playground implementation in our public GitHub repo (see \textit{Data and Code Availability} below). 
\end{itemize}

The remainder of the paper is organized as follows: In Section~\ref{sec:RealWorld}, we provide the real-world context by describing the train rescheduling process at the Swiss Federal Railways (SBB). We describe scope reduction in detail in Section~\ref{sec:researchapproach}. Sections~\ref{sec:ExperimentPipeline} and~\ref{sec:Results} are respectively dedicated to describing the synthetic environment and the computational results. Finally, Section~\ref{sec:Conclusion} concludes.
\section{Real-world Context: Train Rescheduling at the Swiss Federal Railways} \label{sec:RealWorld}

The goal of this section is to describe traffic management in practice and use Swiss Federal Railways (SBB) as an illustration.
Switzerland has one of the world's densest railway networks with both freight and passenger trains running on the same infrastructure. More than 1.2 million people use trains on a daily basis \citep{rcsbrochure}. SBB hence faces challenging rescheduling problems and we use this case to motivate the relevance of our work.

\subsection{Planning Horizon and Control Loop}

\begin{figure}[hbtp]
	\centering
  \includegraphics[width=0.6\textwidth]{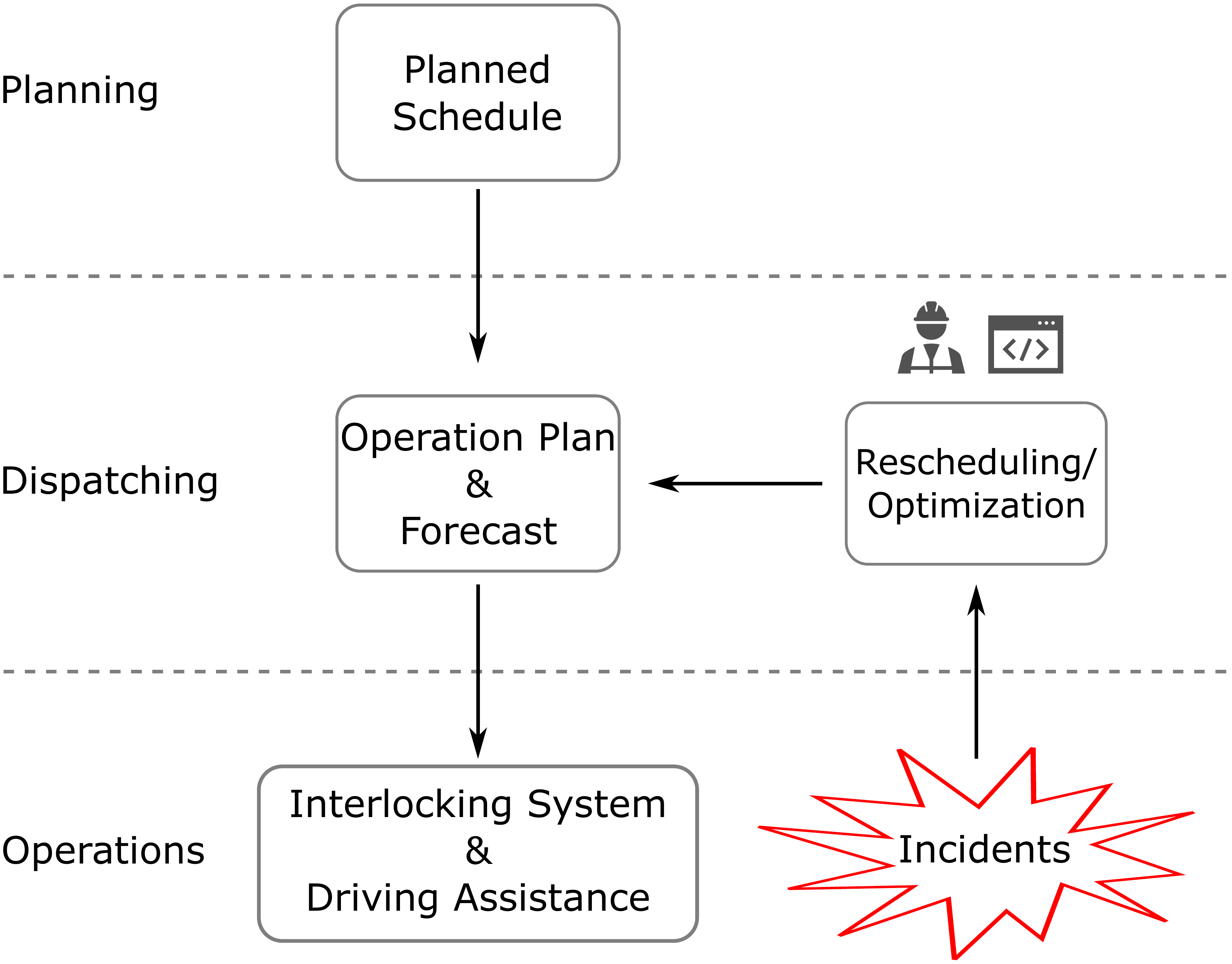}
	\caption{Simplified schema of the rescheduling process during operations. 
	}
	\label{fig:introduction_operations}
\end{figure}

Planning and rescheduling is done at multiple time horizons with their different organizational responsibilities. We depict a simplified view of a rescheduling control loop in Figure~\ref{fig:introduction_operations} \citep[adapted from][]{rcsbrochure,rcswhitepaper}. As we further detail in the following when describing the figure, there are three main areas involved: planning, dispatching and operations. 

The \emph{planned schedule} is the customer interface of the schedule, i.e. the published schedule giving a table of train departures and arrivals that can be used or booked by passengers or train paths that can be reserved by railway undertakings for transportation of goods.

\emph{Dispatchers} have access to information about the current state, and an automated forecast and they make decisions at the level of the operational plan:
\begin{description}
\item[Operational plan] is the railway infrastructure managers internal view of mapping railway services to infrastructure at the microscopic level in the sense of interlocking systems. 
Over a certain time-horizon ahead of production, the operation plan needs to be conflict-free, ensure mutual exclusion of infrastructure usage and additional safety time buffers. To ensure communicability and to ensure service level agreements, decisions should stay closely in agreement with the planned schedule.
At SBB, most of the interlockings can be controlled remotely in four operations centers. However, as the planning process is still largely manual and software tool support is limited to graphical editing and feasibility checking is not fully implemented in the IT systems, the planned schedule has no fully conflict-free associated operational plan.
Furthermore, there  is no global time horizon for the operational plan, conflicts are therefore resolved iteratively.
Nevertheless, in the experiments in this paper, we will assume a conflict-free operational plan.
\item[Forecasts.] Given status updates based on incident reports, there is a forecast of differences and conflicts with respect to the operation plan. Note that not every difference needs to be conflicting, as there may be enough buffer times between trains.
\end{description}

The dispatchers' decisions -- the operational plan -- are translated either automatically or by dedicated dispatching team members to the interlocking system and driving assistance systems. The instructions for the interlocking systems can be time-dependent or event-dependent (for instance, give way to train X after train Y has left). 
At SBB, most of the operation plan is automatically compiled into interlocking instructions. Whereas the interlocking system must be fail-safe and comply with highest safefty standards, the the driving system assistance has only limited safety integrity level, and it only issues recommendations to drivers who remain fully responsible for the actions they take \citep{rcsbrochure}.

\emph{Incidents} (disturbances) are reported as status updates that express small deviations from the operational plan stemming from driving decisions (braking, accelerating), from passenger boarding and alighting, signal box failures, etc. In this work, we do not consider disruptions such as temporary removal of infrastructure, cancellation of trains or connections. 

Over time, new elements from the planned schedule are introduced to the operation plan (rolling time horizon) and the planned schedule also imposes constraints on re-scheduling. For instance, passenger trains must not depart earlier than communicated to customers and the service intention may define further hard or soft constraints such as specific penalties for delay or train dropping.

\emph{Rescheduling/optimization} is the adaptation of the operational plan. It is done both after the update of the forecast based on incidents in the infrastructure, and the unrolling of the time horizon by including future planned schedule into the operational plan. As long as no other trains are affected, the forecast can be directly integrated into the operational plan. The time available to produce a conflict-free solution is typically restricted to a few seconds.

There exist advanced traffic management systems for efficiently re-computing
a traffic forecast for a given time horizon, allowing dispatchers to detect potential conflicts ahead of time. However, the predicted conflicts mostly have to be resolved by humans by explicitly deciding on re-ordering or re-routing of trains based on their experience \citep{VanThiEtAl18}.

\subsection{Hierarchical Organization into Areas of Responsibility}

\begin{figure}[hbtp]
	\centering
  \includegraphics[width=0.6\textwidth]{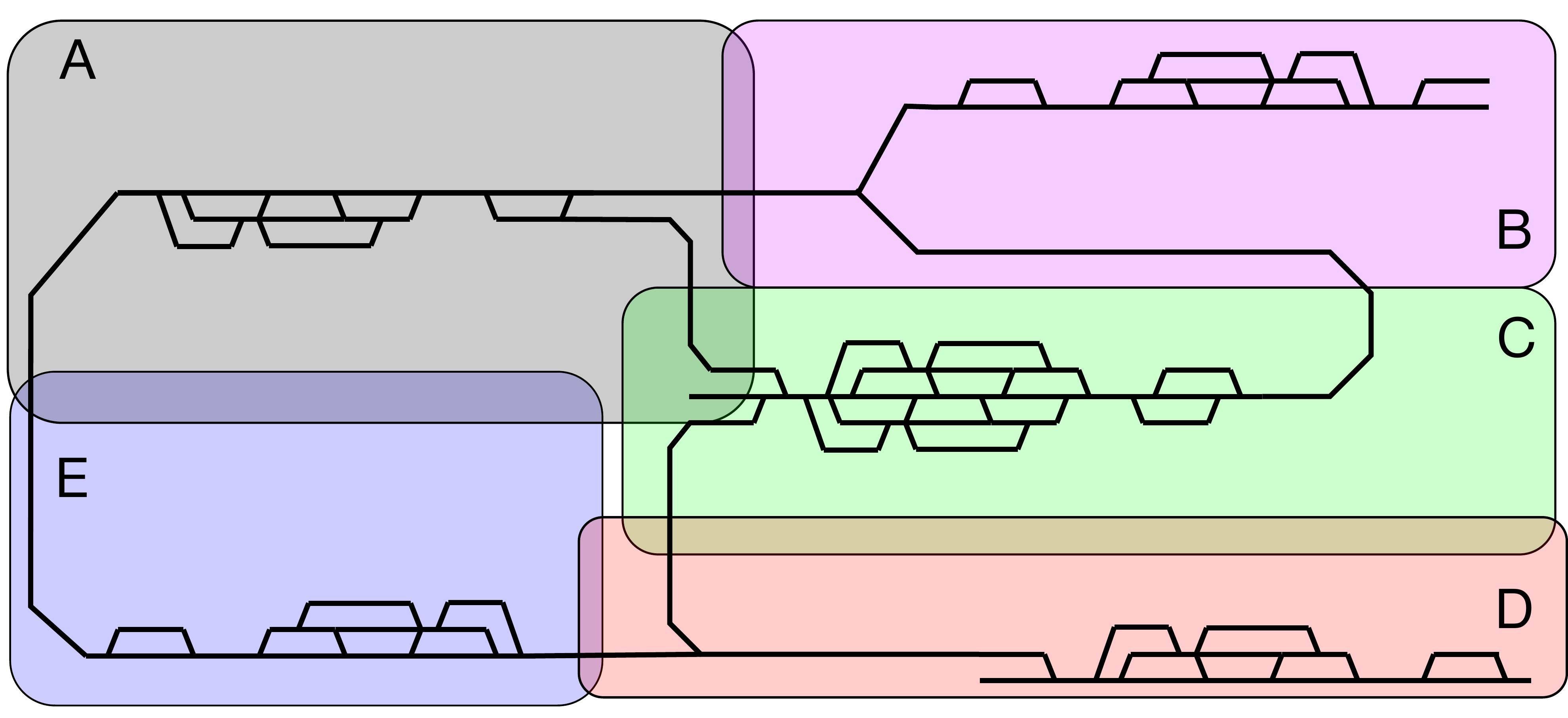}
	\caption{Illustration of geographical decomposition}
	\label{fig:geographical_decomposition}
\end{figure}

Due to the complexity of a railway network, the re-scheduling task is decomposed into smaller geographical areas within which human dispatchers optimize traffic flow. This leads to a multi-level hierarchy of operation centers. In turn, it vastly reduces the complexity within each area and allow human dispatchers to resolve conflicts locally. We show an example Figure~\ref{fig:geographical_decomposition} with five areas (A--E). We intentionally make the areas overlap to illustrate the fact that there is coordination between the areas. Such inter-area coordination is typically done by informal means of communication such as telephone and supported by automated traffic forecasting software. Moreover, dispatchers have access to complete information and can take actions in other regions in agreement with the respective responsible dispatchers.

\subsection{Deployment of Algorithmic Online Re-Scheduling}

At SBB, only very limited areas are fully controlled algorithmically \citep{caimietal2012,Fuchsberger12algorithmsfor,rcsbrochure}. These are small zones around hubs and in tunnels \citep[so-called condensation zones in][]{caimi2009} and current implementations can handle up to one million train re-ordering conflicts per day which are fully resolved by the rail control system \citep{rcsbrochure}.
The goal of such condensation zones is to exploit the hub infrastructure such as approaching areas ahead of major stations or tunnels at capacity. This means, for instance, that trains should enter and pass through the critical sections at full speed. 

Railway networks run on expensive infrastructure and thus an optimal utilization of the given resources is crucial to keep transportation costs low. 
Therefore, automatic control is expected to be deployed to more areas and those areas are expected to grow; thereby, time reserves for catching up and time budget for taking compensation actions will go down,
and coordination between different condensation areas will be required. Humans are not expected to be able to cope with the increased load of decisions to take, and for larger areas, it is currently not possible to find feasible re-scheduling solutions in real time.
The solution space of the train rescheduling problem grows rapidly with important dimensions of the problem, e.g., the number of train routes and the length of the time horizon. Simply deploying the current optimization techniques to a higher number, and larger areas, is not feasible without a scalable strategy for the compensation/coordination work. \cite{VanThiEtAl18} highlight two related limitations in the state of the art. First, exact methods based on decomposition constitute a promising direction but existing approaches are limited in their applicability to large networks. Second, relying on macroscopic approaches may be insufficient since such solutions do not necessarily lead to a feasible microscopic solution.

Based on the above description of the real-world context and the limitations of the state of the art, we believe it is essential that the coordination effort that human train dispatchers make so greatly through communication and experience with an eye on the global system state, can be transported into algorithmic means. Experienced dispatchers have an intuition about which conflicts to pay attention to, while keeping the overall negative impacts under control. We therefore believe that ML could encode the knowledge of defining a restricted optimization scope in both space and time that results in a reduced problem that is computationally tractable with existing solution approaches.

\section{Scope Restriction: Identifying the Core Problem}\label{sec:researchapproach}

In this work, we propose ideas that may lead to a solution approach to overcome the limited spatial restrictions of automated algorithmic re-scheduling by extracting the core problem from the original full problem formulation. We share this preliminary study with the hope that it motivates and encourages further research into the area of real-time problem reduction. In this section, we describe the idea in more detail, and then we outline several different ways to restrict the problem scope. These are used in Section~\ref{sec:Results} to assess the potential of scope restriction.

\subsection{Central Idea}\label{subsec:coreidea}

We describe the idea in a simplified setting. Referring to Figure~\ref{fig:introduction_operations}, we consider an operational plan and a single incident (malfunction) which can cause several conflicts, and we aim to re-schedule such that we have a conflict-free plan again that stays close to the (published) planned schedule. We hint at a fully automatic re-scheduling loop as follows: 
\begin{enumerate}
    \item Incident is found by updating forecast according to real time data.
    \item Core problem is extracted from operational plan and forecast.
    \item Optimizer solves core problem and generates new operational plan.
    \item Repeat (go to~1).
\end{enumerate}

We assume that any re-scheduling problem that arises within the full railway system can be reduced to a much smaller re-scheduling problem defining which trains or even only parts of trainruns might need to be re-scheduled, and we refer to this as the \emph{core problem}. In contrast to the standard approach, we do not follow a spatial decomposition according to network properties but rather investigate the possibility of predicting the relevant core problem given a schedule and a malfunction.

In other words, we assume that it is possible to predict the scope of influence from a given disturbance and that a feasible solution can be found within that scope while keeping the operational schedule fixed outside. The scope consists of: (i) passing times of trains and (ii) routes of trains. It hence spans both space and time. Compared to spatial decomposition, it is closer to a train-based decomposition introduced by \cite{LuanEtAl18}. 

\begin{figure}[hbtp]
	\centering
  \includegraphics[width=0.8\textwidth]{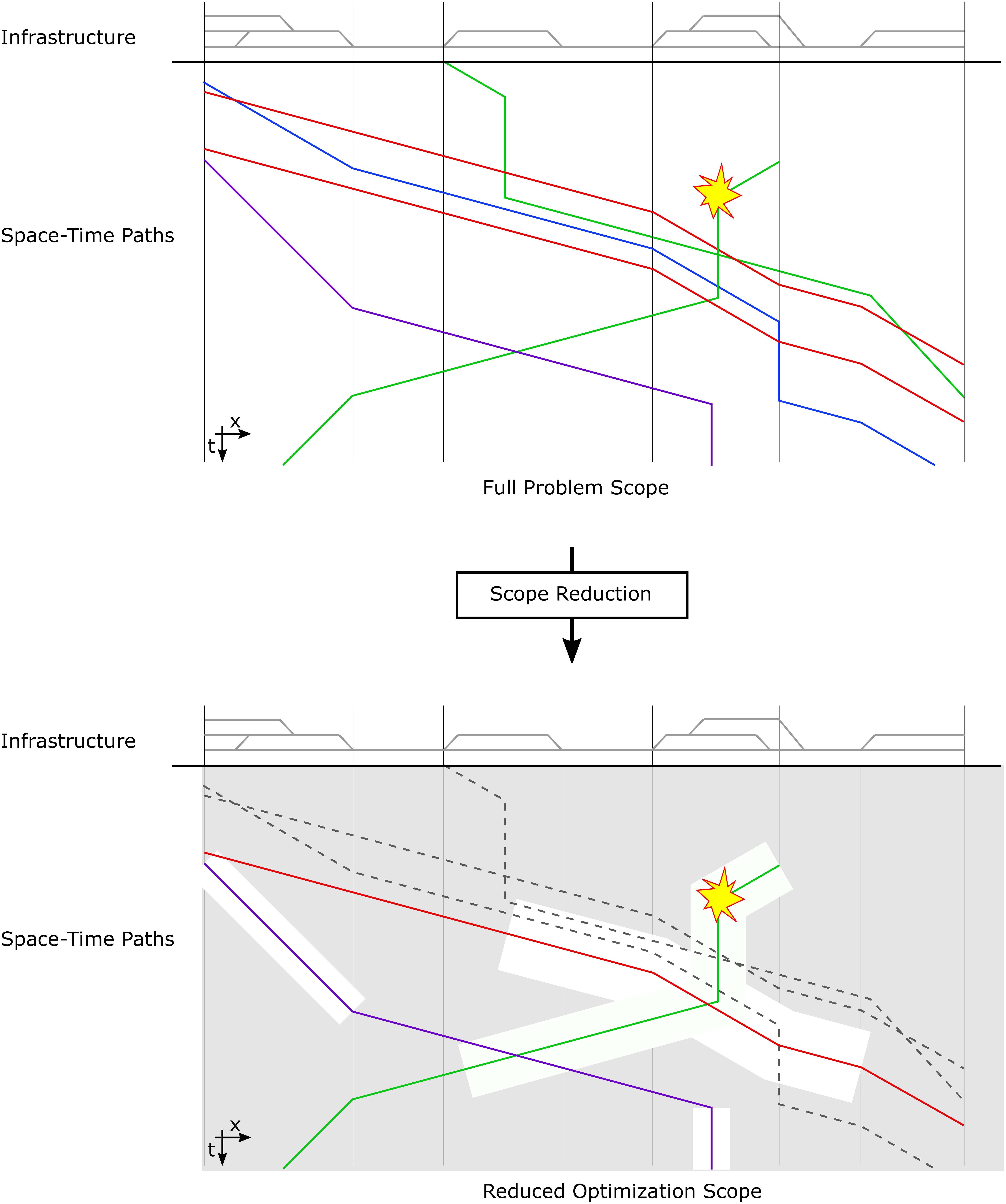}
	\caption{Illustration of scope restriction. The aim is to vastly restrict the search space of an optimizing algorithm, condensing the initial full problem description to the core problem in the time-space parameters. }
	\label{fig:introduction_time_space}
\end{figure}

The main idea of problem scope restriction is shown in Figure~\ref{fig:introduction_time_space}, where a \emph{scoper} predicts which parameters in time and space are relevant to solve the core problem caused by the malfunction. If the scoper can significantly reduce the problem size while allowing for a high-quality feasible solution, this should speed up the solution process. The scoper should hence identify all trains and departures that could be affected by a given disturbance, i.e., that might need adjustment to include the disturbance in the updated operational plan. We note that this problem is related to delay propagation which is a problem that has been extensively studied in the literature \citep[see, e.g.,][]{Goverde10}. 

\subsection{Scopers}\label{subec:online_offline}\label{subsec:scopers}
We now turn our attention to the scopers. As already stated, we believe that an adequately trained machine learning algorithm could constitute a strong scoper. This is left for future research. Instead, in this paper, we focus on assessing the potential of the idea. We therefore define several different scopers that can give useful performance benchmarks.
In this context, we distinguish between two classes of scopers:
\begin{description}
\item[online] The scopers have information about the current state of the trains / network (they do not have access to a full re-schedule solution). This corresponds to the setting that is of interest in practice.

\item[offline] The scopers have information about the current state of the trains / network as well as a full re-schedule solution. These scopers give us useful performance benchmarks.
\end{description}

\begin{figure}[hbtp]
	\centering
  \includegraphics[width=0.4\textwidth]{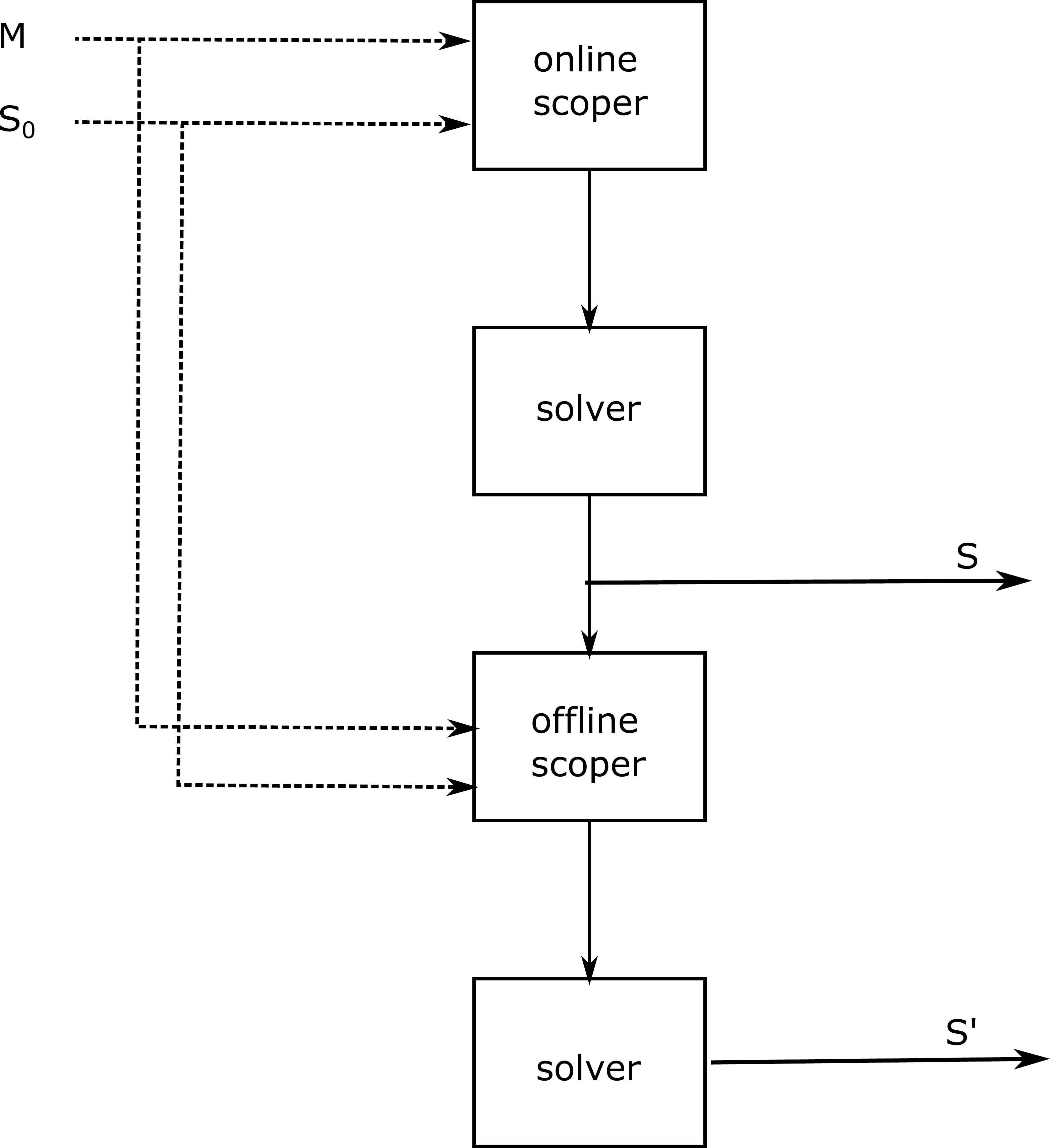}
	\caption{Online versus offline scopers.}
	\label{fig:online_offline}
\end{figure}

We further illustrate the difference between online and offline scopers in Figure~\ref{fig:online_offline}. It shows 
that both online and offline scopers use information about malfunctions ($M$) (state of the trains) and the original schedule ($S_0$). Based on this information, the online scoper defines a core problem that is given to a solver which produces a re-schedule $S$. The offline scoper can use the information of this re-schedule and produce another solution ($S^\prime$). We note that an online scoper that does not restrict the scope (does nothing) is a special case of this scheme.

If we have access to a full solution $S$ of the re-scheduling problem, we can design non-trivial offline scopers that allow for a large speed-ups. Non-trivial means that the output of the offline scoper does not only trivially contain $S$, but a different $S^\prime$ is also a solution. Whereas this offline setting clearly does not correspond to a real application, it provides useful information to assess the potential of the idea. The main question of interest is to assess if it is possible to design an online scoper from data which achieves significant speed-ups with no, or limited impact on solution quality. In addition to the scoper, the answers to these questions depend on the problem instances and the optimization solver. For the experimental results in this preliminary study, we use one specific solver and we mostly focus on analyzing the potential for speed-up using offline scopers.

In the following we describe the different scopers we use for our numerical experiments (see Table~\ref{tab:scopers} for a synopsis). 
More precisely, the online scopers are:
\begin{description}
\item[online\_unrestricted] this scoper does not restrict the re-scheduling problem, it is the ``full'' re-scheduling problem with ``empty'' scope restriction, supposed to give a reference of the hardness of the re-scheduling problem (see $S$ going into \emph{offline scoper} in Figure~\ref{fig:online_offline}); it thus gives a trivial lower bound on speed-up with respect to the solver; this excludes unreachable alternative paths after the malfunction
\item[heuristic] this scoper uses a simple delay propagation algorithm along the scheduled paths to predict which trains will be affected by the malfunction, keeping unaffected trains exactly at their path and times, supposed to show a baseline speed-up (lower bound for non-trivial scopers); we expect the solution quality to deterioriate if there are false negatives.
\item[random] this scoper randomly chooses affected trains, giving no re-routing flexibility to the unaffected trains, supposed to show that the problem of predicting which trains will be affected is not trivial: this is a sanity check online scoping: if we predict affected trains randomly, we expect solutions to be worse than the ones found by the other scopes or; furthermore, we have time flexibility to trains not chosen since for example if a train is scheduled to pass through the malfunction train during its malfunction and is not opened, there is no solution even if we enlarge time windows on the chosen affected trains.
\end{description}

We now turn our attention to offline scopers used in our experiments:
\begin{description}
\item[upper\_bound] this scoper takes the re-scheduling solution from the online\_unrestricted scoper as its scope. This gives a trivial upper bound on speed-up with respect to the specific solver since it measures the specific solver's overhead.
\item[max\_speedup]  this scoper restricts the  scope to the difference between the initial schedule and the online\_unrestricted solution,
\begin{itemize}
    \item only edges from either schedule or online\_unrestricted are allowed;     
    \item if location and time is the same in schedule and full-reschedule, then we stay at them;
\end{itemize}
this is supposed to give a non-trivial but unrealistic baseline on speed-up (upper-bound for non-trivial scopers); in this case, the solution from online\_unrestricted is contained in the solution space, so we expect the same (or an equivalent solution modulo costs) to be found.
\item[baseline] this scoper gives a lose core problem by opening up the same as online\_unrestricted for changed trains. This gives a first impression of how much speed-up we gain by scoping on a train-by-train basis. We expect this to work reasonably well in sparse infrastructures / schedules, but less so in denser infrastructure/schedules, which do not separate well the effect of a malfunction. This might give us hints as to how well a scoper needs to work (with respect to a solver) in order to achieve desired speed-ups.
\end{description}

\begin{table}
\begin{tabular}{|p{36mm}|p{35mm}|p{33mm}|p{43mm}|}
\hline
\thead{scoper}&\thead{routing flexibility}&\thead{time flexibility}&\thead{role}\\
\hline
\hline
online\_unrestricted (online)  &  all trains  &  all trains  &  starting point for offline scopers, measure for hardness of the problem instance\\ \hline
upper\_bound (offline) &  --  &  --  &  trivial upper-bound on speed-up, measuring solver overhead when solution space has exactly one solution \\ \hline
max\_speedup (offline)  &  routes occurring in either schedule or online\_unrestricted solution  &  waypoints (nodes) with a difference between schedule and online\_unrestricted solution  &  close to trivial upper-bound on speed-up\\ \hline
baseline (offline)  &  full routing flexibility for trains with difference online\_unrestricted and schedule  &  trains with difference in online\_unrestricted and schedule  &  a train-based baseline speed-up with access to full solution, giving a non-trivial, but offline upper-bound on speed-up \\ \hline
heuristic (online)  &  all possible paths for trains predicted by transmission chains propagation scheme  &  all nodes for trains predicted by transmission chains  &  lower bound for online scopers \\ \hline
random (online)&  all paths for randomly selected trains (same number as predicted in by heuristic scoper )  &  all nodes from schedule  &  sanity check for baseline, solution quality expected to be worse than in baseline \\ 
\hline
\end{tabular}
\caption{Short-hand synopsis of scopers. Routing flexibility refers to degree of freedom of route choice. 
Time flexibility refers to where along those routes passing times are fixed and where the re-scheduling solution has time flexibility. The last column refers to the scopers' role in our line of argument.}
\label{tab:scopers}
\end{table}

\section{Synthetic Environment and Instance Generation} \label{sec:ExperimentPipeline}

In Section~\ref{sec:Results}, we report results for experiments based on a simplified railway simulation and simple heuristic problem scope restriction algorithms. Our playground implementation described in this section is an abstraction of real railway traffic. 
More precisely, we discuss a synthetic environment based on the Flatland toolbox \citep{mohanty2020flatlandrl} that we use to generate problem instances. This environment is available to other researchers and practitioners. We provide a high-level description here and refer to Appendix~\ref{appendix:implementation} for implementation details. Readers can reach out for more details.
In brief, a given instance of the rescheduling problem corresponds to one train malfunction occurring for a given infrastructure and schedule. We therefore start by describing the infrastructure and the schedule generation followed by the malfunction simulation. Each of these aspects are parameterized to allow a generation of problem instances featuring different characteristics. 

More details on the instance generation can be found in Appendix~\ref{appendix:implementation} and parameters in Appendix~\ref{subsec:experiment_parameters}.

\subsection{Infrastructure}
\begin{figure}[hbtp]
	\centering
  \includegraphics[width=0.8\textwidth]{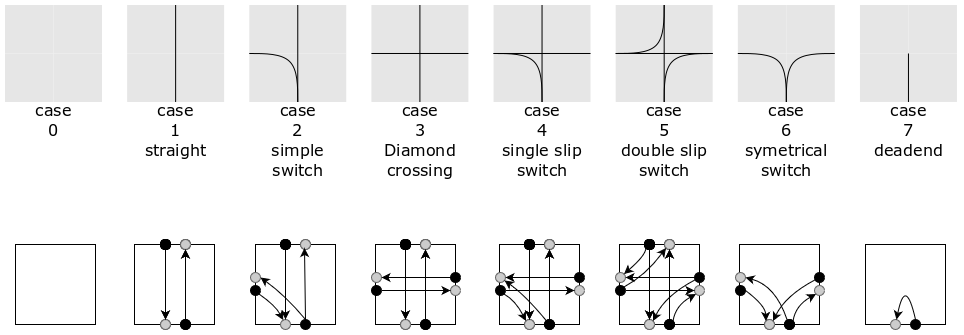}
  \caption{Flatland basic railway elements. The bottom row shows the corresponding translation into a directed double-point graph structure.}
	\label{fig:H1_railway_elements}
\end{figure}

\begin{figure}[hbtp]
	\centering
  \includegraphics[width=0.8\textwidth]{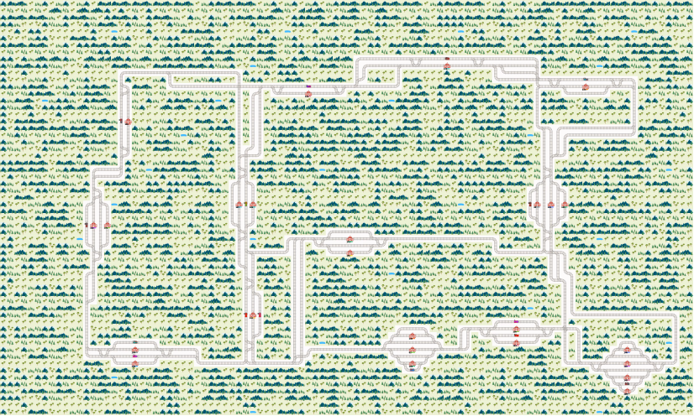}
  \caption{Example infrastructure synthesized by the Flatland generator. It is composed of multiple cities, with randomised number of parallel tracks in cities and between cities.}
	\label{fig:flatland_map}
\end{figure}

We use the Flatland toolbox \citep{mohanty2020flatlandrl} to generate a 2D grid world infrastructure. This world is built out of eight basic railway elements by assigning such an element to each position in the grid. We show the eight elements in the top row of Figure~\ref{fig:H1_railway_elements} each with a corresponding translation into a directed graph in the bottom row. The railway elements define the possible train movements from one position to the four neighbouring ones. From these basic railway elements, it is possible to generate cities 
which in turn are connected together to form a closed infrastructure system as shown in Figure~\ref{fig:flatland_map}. On each side of the city, there are parallel tracks that allow entries and exits, we refer to them as ports in the figure. This structure allows to generate track configurations automatically. More precisely, a fixed number of cities with certain internal track configuration are placed at random on the grid (uniformly distributed). Each city is connected by tracks (given a maximum number of parallel tracks between each city) according to a nearest neighbour algorithm. 

In summary, different infrastructures can be generated by varying the dimension of the grid, the number of cities and the track configuration within and between each city. In turn, these settings can impact the solution space of the re-scheduling problem. For example, parallel tracks increase the re-routing options.

\subsection{Schedule}
For a given specific infrastructure, we generate a train schedule. We take a set of trains $\mathcal{A}$ as given, where each train $a\in \mathcal{A}$ has a single origin and destination in the grid and no intermediate stops. The schedule over time period $t=1,\ldots,T$ defines, for each train, the origin departure time and passing time at each position in the grid (hence also the path in the grid). It is generated such that the sum of train times of all trains is minimized within a given chosen upper bound $T$. We make a number of assumptions: (i) as provided in the Flatland toolbox, each train has a specific but constant speed (ii) there are no time reserves in the schedule, meaning that the trains cannot reduce delay during their journey; (iii) trains appear in the grid when their schedule starts and disappear after reaching their target.

In reality, the schedules are, of course, more complex. For example, there are typically two different schedules: one published and one operational where trains must not depart earlier than published. Moreover, real schedules need to consider a number of aspects that we ignore, for example, respecting fleet management objectives and constraints.

Different schedules can be generated for each given infrastructure by varying the number of trains, their origin-destination pairs and the upper bound $T$.

\subsection{Malfunction Generation}

\begin{figure}[hbtp]
	\centering
  \includegraphics[width=0.6\textwidth]{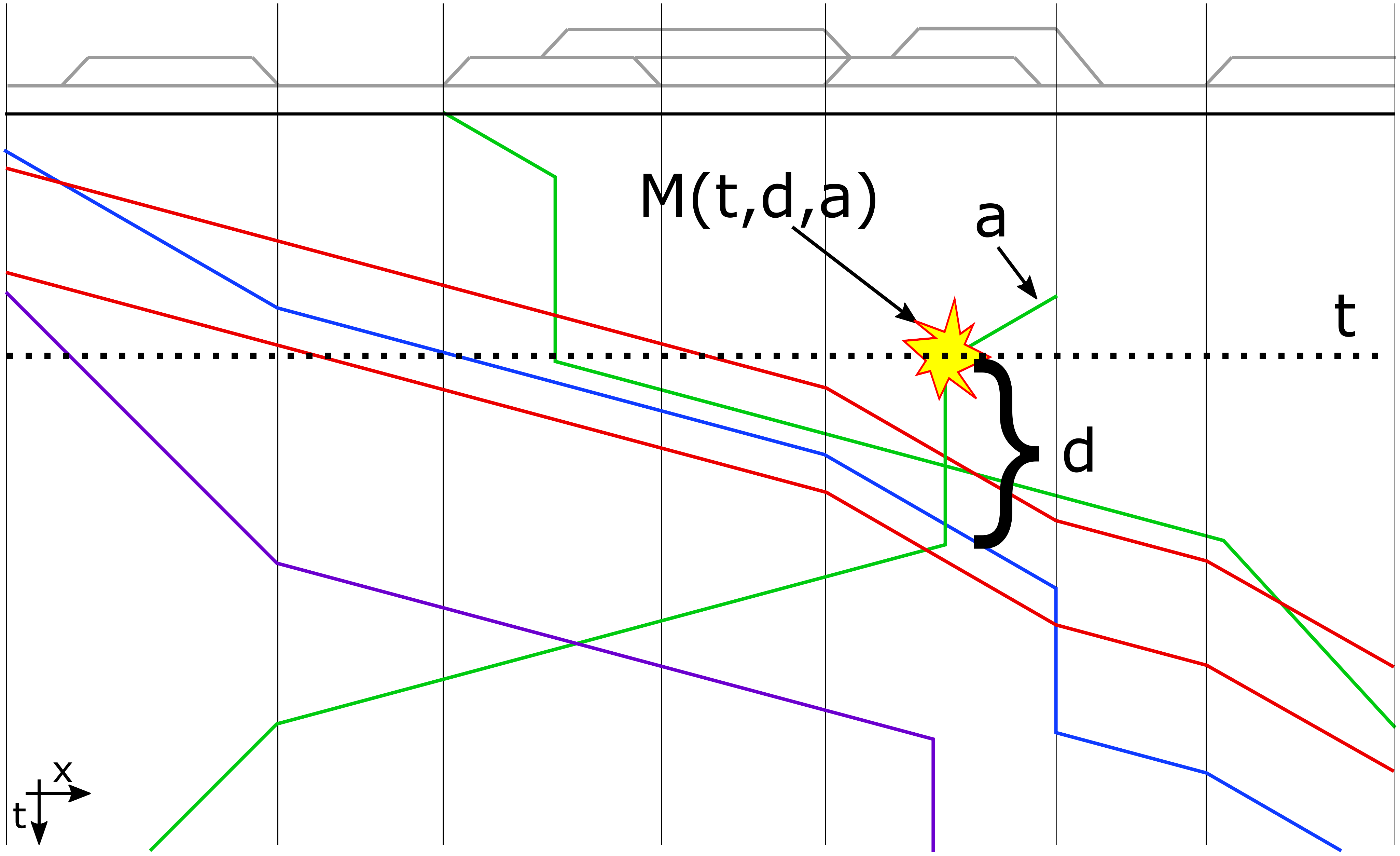}
	\caption{Malfunction generation. This figure illustrates a simple operations plan with 6 trains on a small infrastructure. A single malfunction $M(t,d,a)$ occurs for train $a$ at time $t$ with a total duration of $d$ time steps. }
	\label{fig:introduction_no_loop}
\end{figure}

For a given infrastructure and schedule, we generate different problem instances by simulating malfunctions. In each instance there is one single malfunction $M=(t,d,a)$ and it is generated as follows: We draw a train $a \in \mathcal{A}$ at random, we then draw a malfunction time $t$ at random from $T$. 
The malfunction duration $d$ can be either fixed or drawn at random. 
An example is shown in Figure~\ref{fig:introduction_no_loop} where time is shown on the y-axis (from top to bottom) and space on the x-axis. 

We make two strong simplifying assumptions. First, in reality, multiple malfunctions can occur during the time period $T$ of interest, we assume that only one occurs. This is a strong assumption only if the time period is long or the area of interest is very large. 
Second, at the time of a malfunction, its duration is typically not known in reality while we assume that it is. These assumptions can be relaxed in the synthetic environment but would require different rescheduling models than the ones we consider in this exploratory work.

\section{Computational Results}\label{sec:Results}

In this section, we present the results from our numerical simulations.
The aim of the experiments is to verify the existence of restricted (offline) scopers that allow for a large speed-up in solution time, as well as to investigate the speed-up of non-trivial online scopers. 
The simulations use different scopers for a large set of different schedule and infrastructure configurations. We outline the experimental design in the following section. We then analyze the results first focusing on speed-up (Section~\ref{sec:results_speed_up}), and second, on solution quality (Section~\ref{sec:results_solution_quality}).

\subsection{Experiment Design}

The computational difficulty of any re-scheduling problem depends on the interplay of the infrastructure configuration, the planned railway schedule as well as the location and time of the disturbance. 
As it is hard to characterize what makes an instance hard to solve, we rate the difficulty the instances according to average required computation time to solve the full re-scheduling instance (online unrestricted scoper). 

Using an hierarchical setup for problem instance generation, as shown in Table~\ref{tab:hierarchical_experiment_design}, we produce a diverse set of problem instances with varying levels of difficulty (as measured by computation time). The full parameter ranges are shown in Appendix~\ref{subsec:experiment_parameters} and the distribution of instance difficulty can be found in Appendix~\ref{appendix:exp_results}.

When reporting the results in the following sections,  we order instances according to the  relative difficulty of the re-scheduling instances.
While schedule and infrastructure remain mostly static in real world railway networks, we chose to generate different infrastructure configurations to get a wider distribution of problem instances.
However, it turned out to be difficult to identify the relevant drivers for the instance difficulty. The resulting distribution is biased towards low computing time (more than $60\%$ of experiments are in the shortest bin with computation times of less than $38$ seconds).

All results were obtained by using the publicly available model of \cite{DBLP:journals/corr/abs-2003-08598,
DBLP:conf/lpnmr/AbelsJOSTW19}\footnote{\url{https://github.com/potassco/train-scheduling-with-clingo-dl}} using an Answer Set Programming (ASP) solver, without specific tuning of solver parameters. 
The solver objective penalizes delay and minimizes the number of changes in the rescheduling solution\footnote{More precisely, the objective for re-scheduling is a weighted sum of the delay at the target and a penalty for the number of nodes in the graph the re-scheduled path leaves to the scheduled path. The last term should help to avoid ``flickering'', the
re-scheduled should not deviate from the scheduled path ``without need''.}.
We used the same solver parameters for all instances, and the computations were run in a cloud environment.

\begin{table}[ht]
    \centering
    \begin{tabular}{p{1cm}|p{2cm}|p{11cm}}
        \textbf{Level} & \textbf{Element} & \textbf{Variables} \\ \hline
       0 & infrastructure  &  Dimension of the grid as well as density of train tracks and number of trains.\\ \hline
        1 & schedule & Schedule times and paths for all trains.\\ \hline
        2 & malfunction & Malfunctioning train and malfunction onset.\\ \hline
        3 & solver runs & Each configuration was run multiple times with different seeds 
    \end{tabular}
    \caption{Hierarchical experiment design. For any level, all subsequent levels generate unique problem instance configurations.} 
    \label{tab:hierarchical_experiment_design}
\end{table}


\subsection{Speed-Up}\label{sec:results_speed_up}

\begin{figure}[hbt]
            \includegraphics[width=\textwidth]{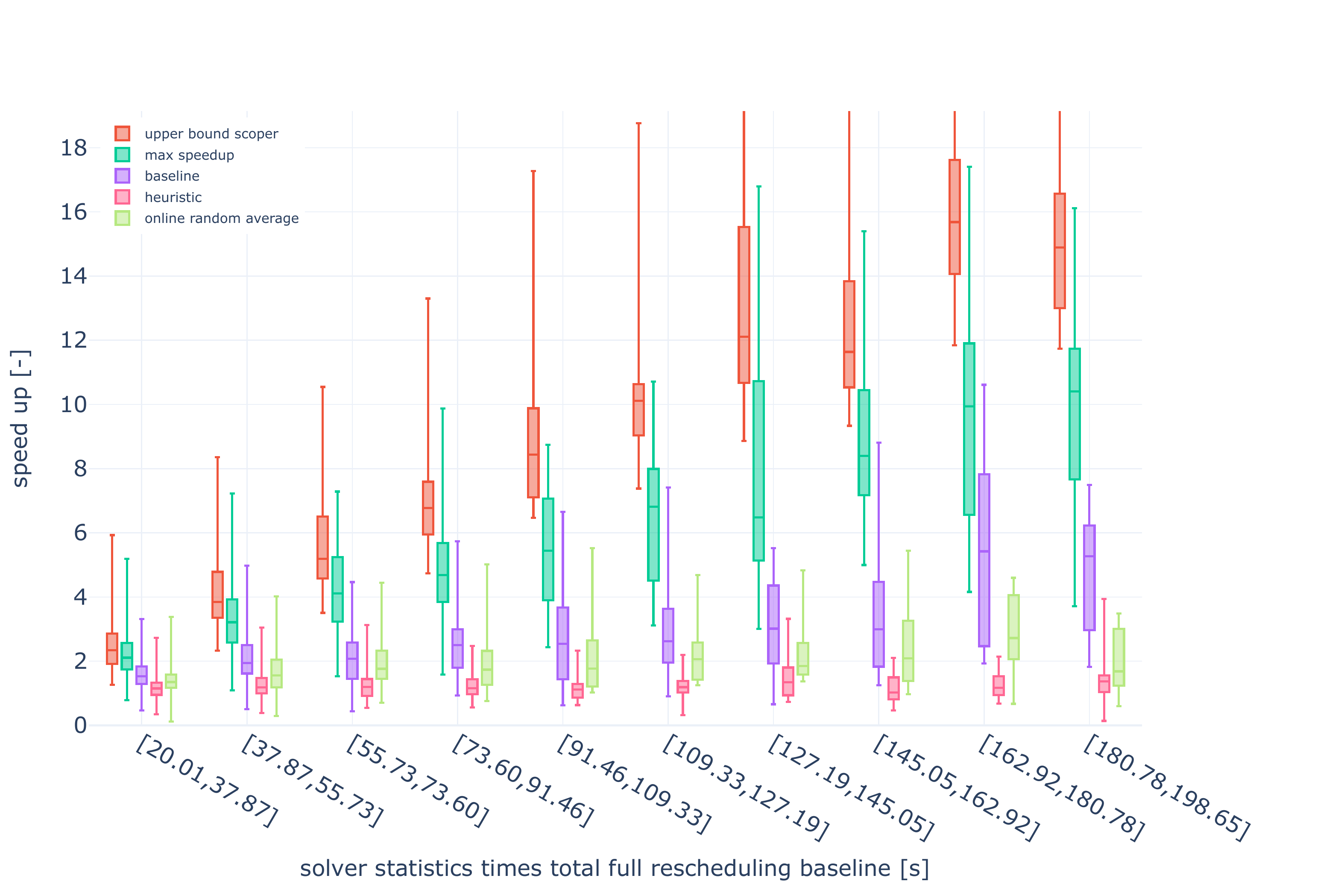}
         \label{fig:computationtimes_speed_up_total}
	\caption{Speed-up, defined as the the factor between the computation time needed to solve the full scope and the restricted scope (per scoper), where a higher number means more improvement.
	}
	\label{fig:speed_up}
\end{figure}

In this section, we analyze the speed-up, that is the reduction in computation time needed by the solver\footnote{We measure elapsed time, including the grounding phase of the ASP solver, without optimizing the encoding to speed up neither the grounding nor the solving phase. The absolute computation times and relative speed-ups of the scope restriction may be different for other ASP encodings and other OR solvers. We only investigate the speed-up in this one setting.} after the scope restriction compared to the one without scope restriction. The results are shown in Figure~\ref{fig:speed_up} where we only analyze experiments in the range of interest for real-world applications of $[20s,200s]$ computation time for the full problem instance. The results are binned in $10$ equidistant bins for the analysis. We omit a comparison of absolute computational time for rescheduling  and refer the interested reader to the Appendix~\ref{appendix:exp_results}.

The \textit{upper\_bound scoper} is a sanity check for our approach. The results illustrate the maximum speed-up that can be achieved with the ASP-Solver. The speed-up is limited by solver specific pre-processing. This check is achieved by locking all variables to the final solution and letting the ASP-Solver build and solve the problem instance without any degrees of freedom. The computational cost of this is specific for the used ASP solver and can, of course,  vary for other solvers.

Our estimate for maximum speed-up factors is represented by the \textit{max speedup} scoper, which consistently achieves a factor above $2$ and up to $10$. Given that this scoper makes use of the full problem information to deduce the minimal core problem it constitutes a maximum in our setting.

The \textit{baseline scoper} uses a much looser restriction on the scope. It restricts the scope to only the affected trains, but provides full routing and time flexibility. An accurate online scoper should hence be able to achieve this speed-up.
We observe speed-up factors ranging between $2$ and $5$. These are large enough to have a positive impact on real-world operations.

The \textit{heuristic scoper} is desiged as a first attempt at an online scoper. It uses a simple delay propagation logic to define the core problem. Important in this context is the classical compromise between false positive and false negative rates. On the one hand, if the scoper has a high false positive rate, it will lead to little or no speed up. On the other hand, if the false negative rate is high, this can have a detrimental impact on solution quality which we analyze in the following section. 

Based on the results reported in Figure~\ref{fig:speed_up} we clearly see the impact of a relatively high false negative rate. Our simple heuristic scoper does not achieve a significant speed up.

As a sanity check we also report results for \textit{random average scoper}, which is the average solution time of using $5$ randomly generated \textit{core problems}. This scoper only shows a minor speed-up. 

The results support our hypothesis that sufficient decrease in problem scope can lead to significant speed-up in computation time. These findings are expected in light of the literature on heuristic variable fixing. Furthermore, we see that, at least for the given simplified rail infrastructure and schedules, the core re-scheduling problem is much smaller than the full re-scheduling problem. This even holds true for comparably small railway networks and simple schedules. However, identifying the core problem is not a trivial task.

\subsection{Solution Quality}\label{sec:results_solution_quality}

\begin{figure}[tb]
            \includegraphics[width=\textwidth]{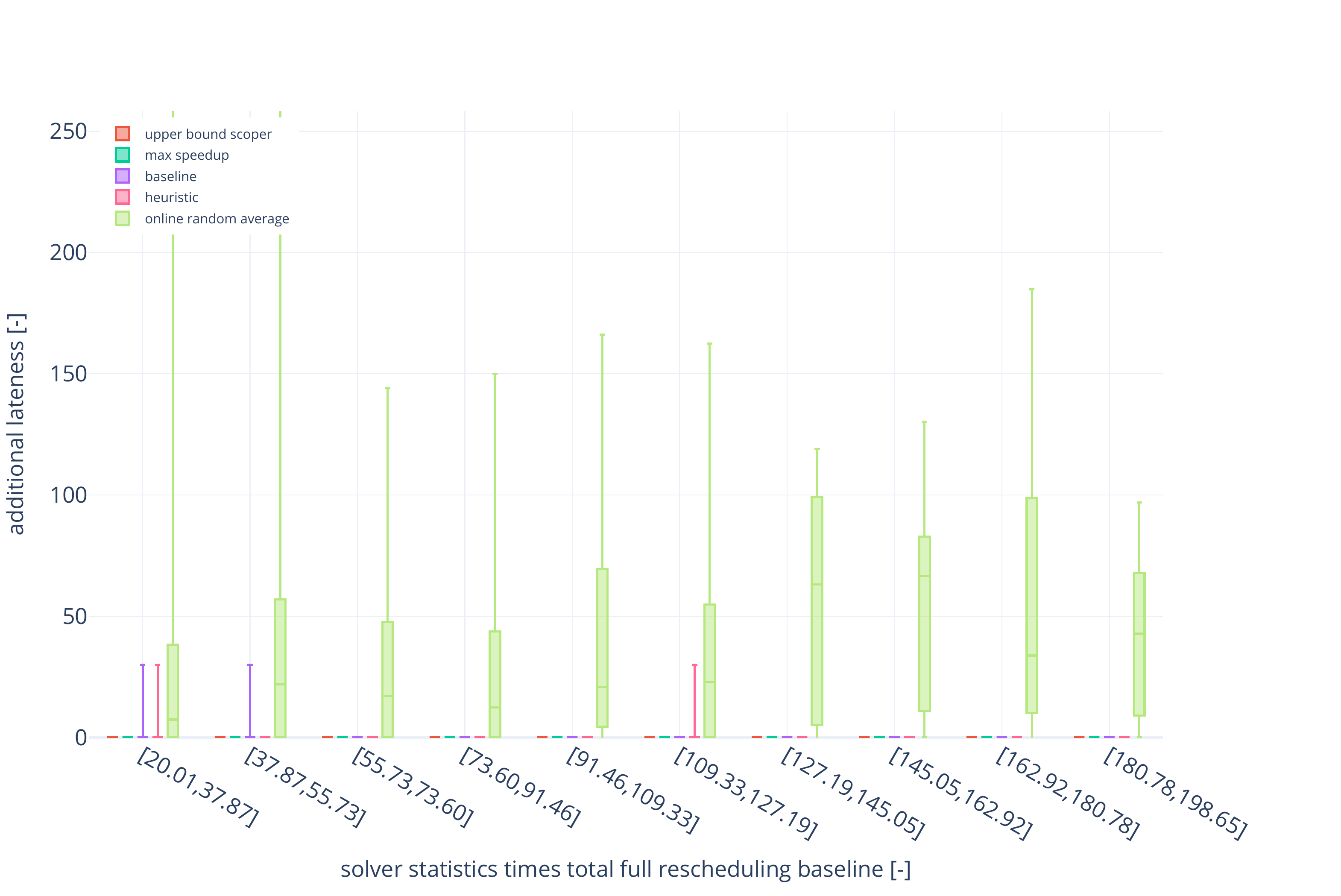}
            \caption{Solution quality. This plot shows the additional lateness that was generated compared to the perfect rescheduling solution. 
            }
	\label{fig:delay}
\end{figure}

In this section we briefly comment on solution quality as defined by the objective function. In our case it consists in minimizing the temporal delay compared to the original schedule. Given our simulation environment, we do not investigate solution quality in more general terms taking into account business rules and industry metrics.

Figure~\ref{fig:delay} 
 shows an aggregate view of lateness (compared to a perfect rescheduling solution) for all instances. We note that all but the random scoper have high solution quality. By design, the baseline and max speed-up scopers yield a solution of the same cost as the optimal solution. The additional lateness in the case of the baseline scoper is compensated by fewer track changes (not shown).  For all but a few instances, our heuristic scoper also yields the same high quality solutions. It can be explained by its low false negative rate. Recall, however, that this scoper did not yield any speed up due to its high false positive rate.
 As expected, the random scoper, achieves a poor solution quality. 
 
 These preliminary results show that, according to expectation, we can achieve a relatively high speed up by identifying trains affected by a malfunction with little negative effect on solution quality (baseline scoper). However, as we already stated, accurately defining the core problem is a non-trivial task.

\section{Conclusion}\label{sec:Conclusion}

Effective rail-time railway traffic management is essential for reliable train services.
In practice, operators rely on a large extent on humans and their expertise to make rescheduling decisions. These human dispatchers heavily rely on direct human communication to coordinate their decisions. On the contrary, algorithmic solutions in industry are only implemented on restricted geographic scopes due to computational budget restrictions. In turn, the geographic decomposition introduces coordination challenges and is therefore difficult to scale.

Motivated by results in the literature and practical challenges, we argued that it should be possible to define a scoper (i.e., an algorithm) that can predict a reduced problem scope. We referred to this as the core problem which restricts the original problem in time and space by identifying train services that need to be rescheduled. The aim is to find high-quality traffic management solutions within a limited computing time budget and without introducing a coordination challenge.

We reported preliminary results from an exploratory study where we used the Flatland simulation environment to generate problem instances. Based on the results we concluded that, as expected, accurately identifying the core problem can lead to important speedups. However, identifying the core problem is a non-trivial task and 
a high-performing scoper needs to avoid false negatives (to ensure feasibility resp. solution quality) while not having too many false positives (to keep speed-up)
when it comes to identifying the trains that are affected by a disturbance.
This was an exploratory study and much research is needed to turn these ideas into an operational approach. How to define an accurate scoper\footnote{At a more theoretical level, how many equivalent core problems are there and how to decide between them?}, how to best interact with a given solver, how to guarantee feasibility, and how to characterize uncertainty around the core problem are examples of four open and essential questions. Furthermore, not all of the specific business rules need to be reflected by the solver cost function -- the solver and its metric might be part of a larger pipeline in an industrial setting.

The main purpose of sharing these preliminary ideas and results are to stimulate more work on this research topic. We believe that adequately defined machine learning algorithms trained on historical and/or simulated data should be able to accurately predict the core problem. The main reason is that infrastructure and train services are relatively stable over time. We hence face a recurrent problem that experienced humans are good at solving. We provide an extensible playground open-source implementation and our benchmarking instances and data in our GitHub repository (see \emph{Data and Code Availability} Section below).

\section*{Data and Code Availability}
\addcontentsline{toc}{section}{Data and Code Availability}
Coda and data are available in our public GitHub repositories\footnote{\url{https://github.com/SchweizerischeBundesbahnen/rsp}}$^,$\footnote{\url{https://github.com/SchweizerischeBundesbahnen/rsp-data}}. Also, the paper with an extended technical appendix can be found in the code repository.

\section*{Author Contributions Statement}

\addcontentsline{toc}{section}{Author Contributions Statement}
Following CRediT (Contributor Roles Taxonomy\footnote{\url{https://credit.niso.org/}}), the authors have contributed as follows.
\begin{description}
    \item[Emma Frejinger] Conceptualization, Writing (Literature Research), Writing (Review and Editing)
    \item[Erik Nygren] Conceptualization, Methodology, Software (Analysis), Writing (Introduction, Computational Results, Review and Editing)
    \item[Christian Eichenberger] Software (Pipeline and Analysis), Conceptualization, Methodology, Writing (First Integral Draft, Review and Editing)
\end{description}

\section*{Acknowledgement}
\addcontentsline{toc}{section}{Acknowledgement}

This research project was only possible due to the generous support of both Mila and CIRRELT which provided access to working places, compute and the opportunity for great scientific exchanges.

We in particular want to thank Andrea Lodi and Yoshua Bengio who helped shape this idea through discussions and feedback.

Furthermore we want to thank the Swiss Federal Railway company for creating this opportunity for an intensive international research collaboration by providing resources as well as access to their core challenges. In particular we would like to thank Dirk Abels, Mathias Becher and J{\"u}rg Balsiger for supporting the research efforts and paving the way for an unobstructed research exchange between SBB, Mila and CIRRELT. Also we would like to thank for the fruitful discussions with colleagues at SBB (PFI, Flux, RTO) and University of Potsdam (Potassco).


\bibliographystyle{plainnat_custom}
\bibliography{biblio}

\appendix

\begin{appendices}

\section{Implementation Overview for the Hierarchical Experiment Setup}\label{appendix:implementation}

We now give an overview of the pipeline from an implementation perspective.

\begin{figure*}[p]
\begin{subfigure}[t]{0.47\textwidth}
\centering
\includegraphics[scale=0.2]{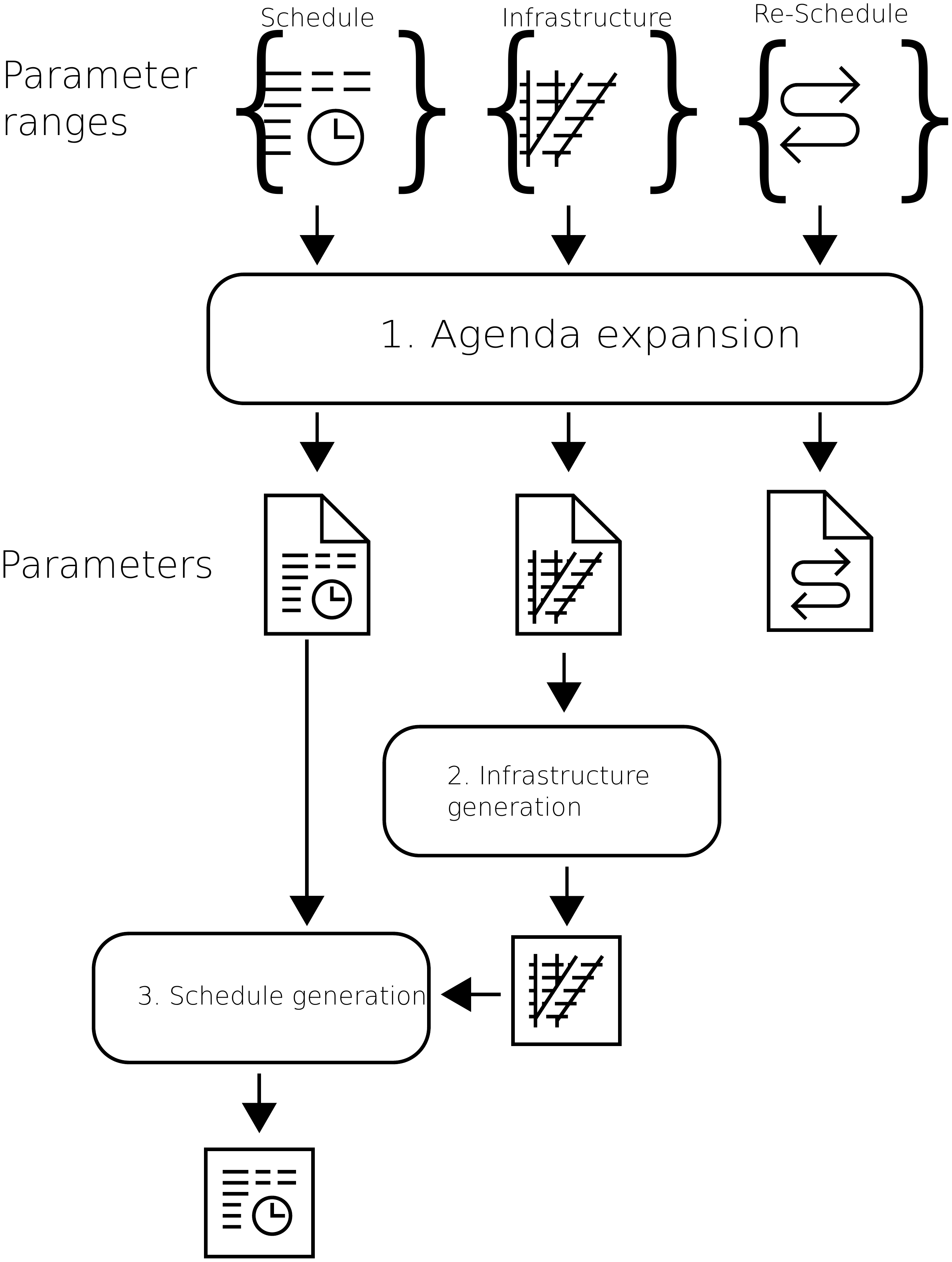}
\caption{Offline preparation for pipeline. The output is a hierarchical directory structure, where each infrastructure can have multiple schedule and each schedule multiple reschedule parameters. These files are generated by expanded initial parameter ranges for schedule, infrastructure and re-scheduling.}
\label{fig:H1_overview}
\end{subfigure}
\hfill
\begin{subfigure}[t]{0.47\textwidth}
\centering
\includegraphics[scale=0.2]{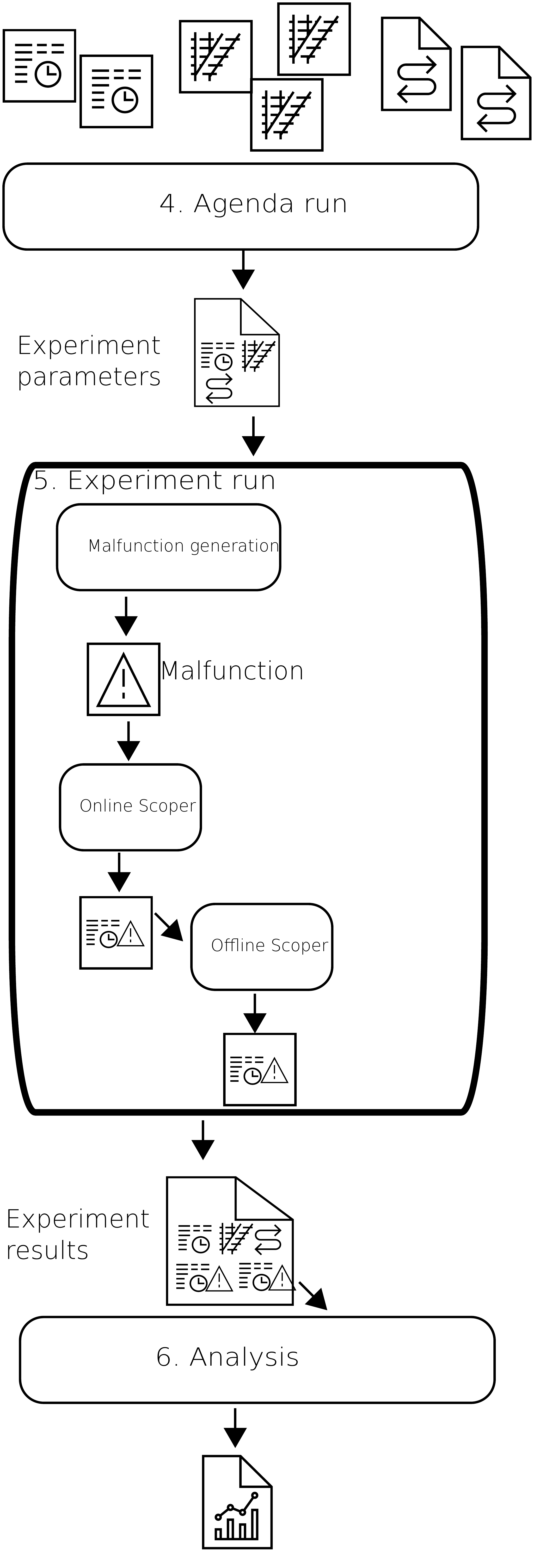}
\caption{Agenda run starts from the offline preparation hierarchy of infrastructures, schedules and re-schedule parameters (see Figure~\ref{fig:H1_overview}) from which a subset is chosen. An agenda run consists of multiple experiment runs. 
Each single experiment run gets experiment parameters as input and runs different scopers on them. The resulting re-schedules and performance measure are stored in Experiment Results. The results from multiple runs can be visualized by the analysis component.
}
\label{fig:experiment_pipeline}
\end{subfigure}
\caption{Rescheduling experiments pipeline: preparation, experiments, analysis. Cornered rectangles represent files or data structures and rounded rectangles represent steps.}
\end{figure*}

Since the schedule generation step takes too long to repeat for every experiment, we will use the same infrastructure and schedule for many malfunctions. In fact, we can pre-generate infrastructures and schedules and then work on variations of the re-scheduling part of the pipeline more efficiently. 
Referring to Figures~\ref{fig:H1_overview} and \ref{fig:experiment_pipeline}, the pipeline decomposes into five top-level stages plus a post-experiment analysis stage:
\begin{description}[%
  before={\setcounter{descriptcount}{0}},%
  ,font=\bfseries\stepcounter{descriptcount}\thedescriptcount~]
\item[Agenda Expansion] The  parameter  ranges  are  given  as  input  and  expanded into infrastructure and schedule parameters to generate infrastructure and schedule.
\item[Infrastructure Generation] This generates railway topologies and places train start and targets in the infra\-structure\footnote{Conceptually, the placement of trains should be separated from infrastructure generation. However, the way Flatland is designed, it is not possible without refactoring to separate the two as the information about cities/stations is only available during infrastructure generation and train placement.}.
\item[Schedule Generation] This generates the exact conflict-free paths and times through the infrastructure for all trains. 
\item[Agenda Run]  From hierarchical directory of generated infrastructures, generated schedules and re-schedule parameters, an agenda of experiments is compiled and experiment are run sequentially.
\item[Experiment Run] Here, the different scopers are run on the same problem instance as we know it from Figure~\ref{fig:online_offline}.
\item[Agenda and Experiment Analysis] This stage aggregates experiment results from multiple experiment results and generates the plots 
in Section~\ref{sec:Results}.
\end{description}

\clearpage

\section{Experiment Parameters}\label{subsec:experiment_parameters}

Here we give more details on the experiments for the results of Section~\ref{sec:Results}.
We use $3264$ experiments from a grid spanned by the parameter ranges below.
We have 12 infrastructures (Level 0), 4 schedules per infrastructure (Level 1), $(50+62+74+86)$ malfunctions (Level 2),  and only have one solver run per experiment (Level 3).\footnote{The exposition of the parameter ranges slightly deviates from the structure in the data repository. Instead of 12 infrastructures with 4 schedules, technically, we have 48 infrastructures with 1 schedule each (using \texttt{flatland\_seed\_value}=$[190,191,4]$) . Note that, due to multi-threading, schedule generation is non-deterministic even when using the same ASP seed.}
The value ranges $[a, b, n]$ have the following meaning: 
\begin{itemize}
    \item if $n=1$, we take only $a$ (ignoring $b$); 
    \item if $n>1$, we take $n$ points from the half-open interval $[a,b)$ of step size $\lfloor \frac{b-a}{n}\rfloor$.
\end{itemize}
In the following tables, we abbreviate $100$ instead for $[100, 100, 1]$
and $[0,48]$ for $[0,48,48]$. Examples: $[8,15,3]$ expands to $8,10,12$, and $[1,2,2]$ expands to $1,1$.

\subsection{Infrastructure Parameters}\label{sec:infraparams}
\begin{tabular}{|p{40mm}|p{15mm}|p{40mm}|p{25mm}|}
\hline
\thead{Parameter} & \thead{Symbol} & \thead{Description} & \thead{Value}\\
\hline
\hline
\texttt{infra\_id} & -- & infrastructure id & $[0,48]$\\\hline
\texttt{width} & $w$ & Number of cells in the width of the environment & 100\\\hline
\texttt{height} & $h$ & Number of cells in the height of the environment & 100 \\\hline
\texttt{flatland\_seed\_value} & -- & Random seed to generate different configurations & [190,190,1] \\\hline
\texttt{max\_num\_cities} & $\left|S\right|$& Maximum number of cities to be places in the environment. Cities are the only places where trains can start or end their journey. Cities consists of parallel track and entry/exit ports. & $[8, 15, 3]$\\\hline
\texttt{grid\_mode} & -- & & \texttt{False}\\\hline
\texttt{max\_rail\_between\_cities} & -- & Maximum number of parallel track at entry/exit ports of the cities & $1$\\\hline
\texttt{max\_rail\_in\_city} & -- & Maximum number of parallel tracks in the city & $2$\\\hline
\texttt{number\_of\_agents} & $\left|\mathcal{A}\right|$& & $[62, 86, 4]$\\\hline
\texttt{speed\_data} & $v: \mathcal{A}\to [0,1]$ & distrbution of speeds among trains & \texttt{\{
            1: 0.25,
            1/2: 0.25,
            1/3: 0.25,
            1/4: 0.25
        \}}\\\hline
\texttt{\mbox{number\_of\_shortest\_paths} \mbox{  \_per\_train}} & $\approx(\mathcal{V}_a,\mathcal{E}_a)$ & We compute shortest paths only once; should be larger than the number in scheduling and re-scheduling.& 10 \\
\hline
\end{tabular}

\subsection{Schedule Parameters}\label{sec:schparams}

\begin{tabular}{|p{40mm}|p{15mm}|p{40mm}|p{15mm}|}
\hline
\thead{Parameter} & \thead{Symbol} & \thead{Description} & \thead{Value}\\
\hline\hline
\texttt{infra\_id} & -- & reference to infrastructure & $[0,48]$ \\\hline
\texttt{schedule\_id} & -- & schedule id& $[0,4]$\\\hline
\texttt{asp\_seed\_value} & -- & Since we use 2 threads, the ASP solver behaves non-deterministically, so the seed value has no effect.& $814$ \\\hline
\texttt{number\_of\_shortest\_paths \_per\_train\_schedule} & $\approx(\mathcal{V}_a,\mathcal{E}_a)$ &  & $1$\\
\hline
\end{tabular}

\subsection{Reschedule Parameters}\label{sec:reschparams}

\begin{tabular}{|p{40mm}|p{15mm}|p{40mm}|p{15mm}|}
\hline
\thead{Parameter} & \thead{Symbol} & \thead{Description} & \thead{Value}\\
\hline
\hline
\texttt{earliest\_malfunction} & $m_{earliest}$ & Used to determine $m_{time\_step}$, the time step of the malfunction. & $30$\\\hline
\texttt{malfunction\_duration} & $m_{duration}$& Malfunction duration.& $50$\\\hline
\texttt{malfunction\_train\_id} & $m_{train}$& Which train is disturbed? &$[0, 86]$\\\hline
\texttt{number\_of\_shortest\_paths \_per\_train} & $\approx(S,L)$& Defines the route restrictions for the trains.& $10$\\\hline
\texttt{max\_window\_size\_ from\_earliest} & $c$& Truncate window sizes to this time window sizes to this maximum. Applies to windows on vertices and not on resources. & $60$\\\hline
\texttt{asp\_seed\_value} & -- & Since we use 2 threads, the ASP solver behaves non-deterministically, so the seed value has no effect. &  99\\\hline
\texttt{weight\_route\_change} & $\rho$ & How many time steps delay is one route change equivalent to? & 30 \\\hline
\texttt{weight\_lateness\_seconds} & $\delta_{penalty}$ & Factor to scale delays. Should only be different than $1$ if \texttt{weight\_route\_change} is equal to 1, and vice-versa. & $1$ \\
\hline
\end{tabular}

\clearpage
\section{Further Experiment Result Details}\label{appendix:exp_results}
\subsection{Re-scheduling Times Distribution}
\begin{figure}[hbtp]
    \includegraphics[width=\textwidth]{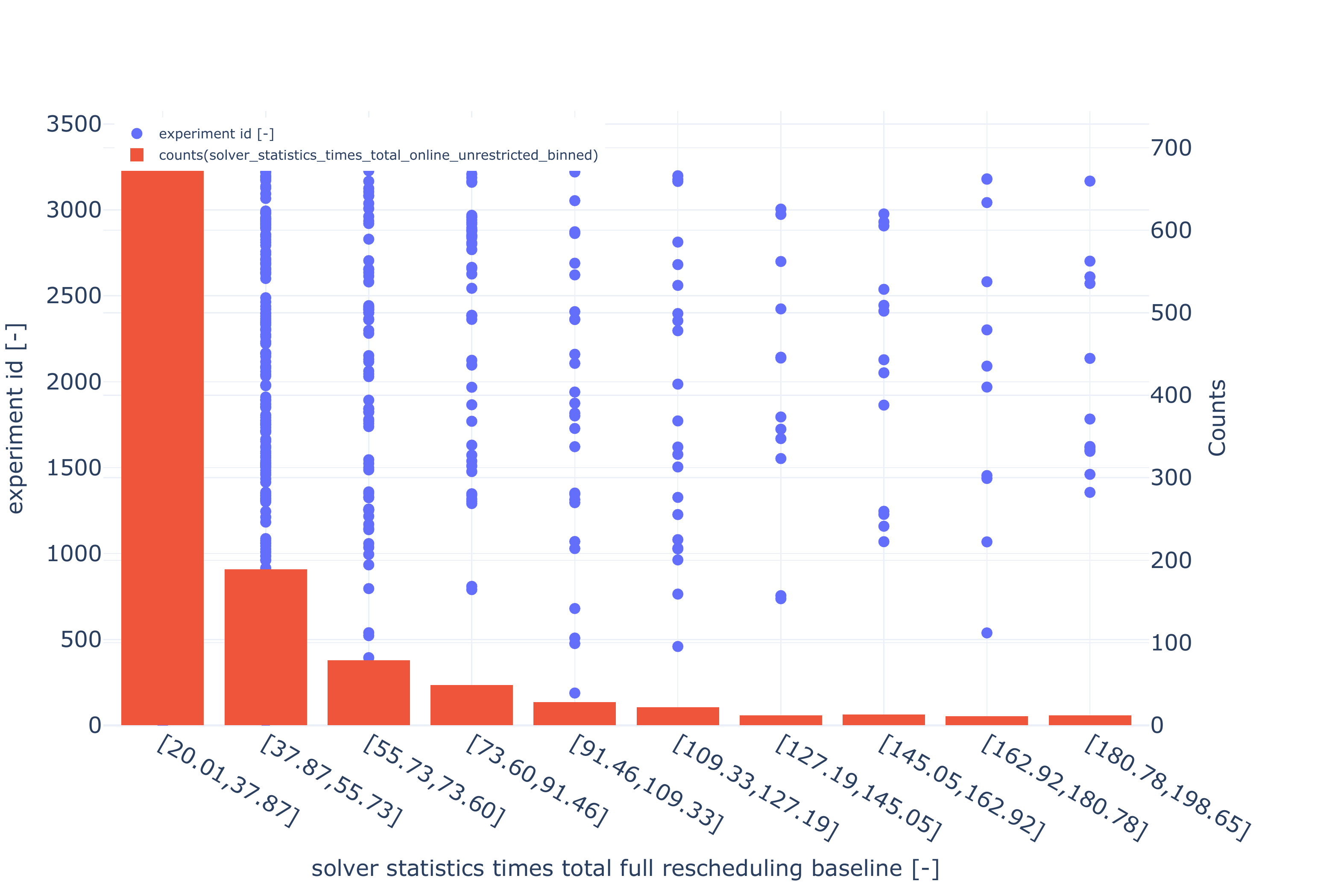}
	\caption{Histogram of experiments in 10 equidistant bins, filtered on data with \texttt{solver\_statistics\_times\_total\_online\_unrestricted} between $20s$ and $200s$. There is no correlation between experiment ID and computation time, meaning that no single factor (e.g. infrastructure, schedule or malfunction) is defining the complexity of a problem instance. The histogram further shows, that most problem instances have low complexity (short solving time) where our approach will have little impact.}
	\label{fig:agenda_histogram}
\end{figure}

Figure~\ref{fig:agenda_histogram} shows that the hardness of re-scheduling is not purely determined by the infrastructure parameters (in particular grid size), but that there is a fat tail of run times depending on the malfunction configuration within the infrastructure.

\begin{figure}[hbtp]
   \includegraphics[width=\textwidth]{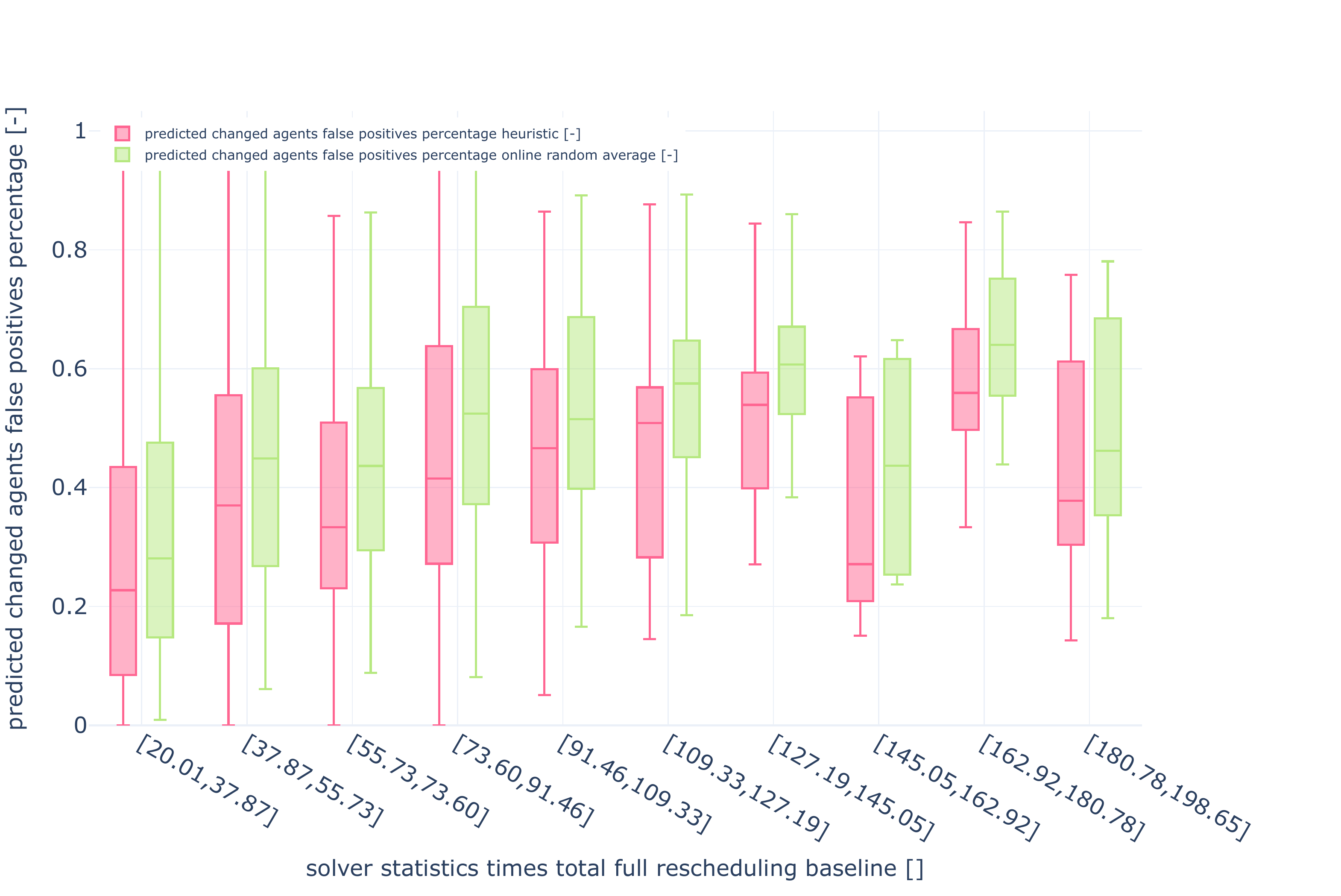}
	\caption{Prediction quality: false positives and false negatives rate. Online transmission chains fully restricted is not shown as it is almost completely identical to online transmission chains route restricted.}
		\label{fig:prediction_quality_false_positive}
\end{figure}

\begin{figure}[hbtp]
   \includegraphics[width=\textwidth]{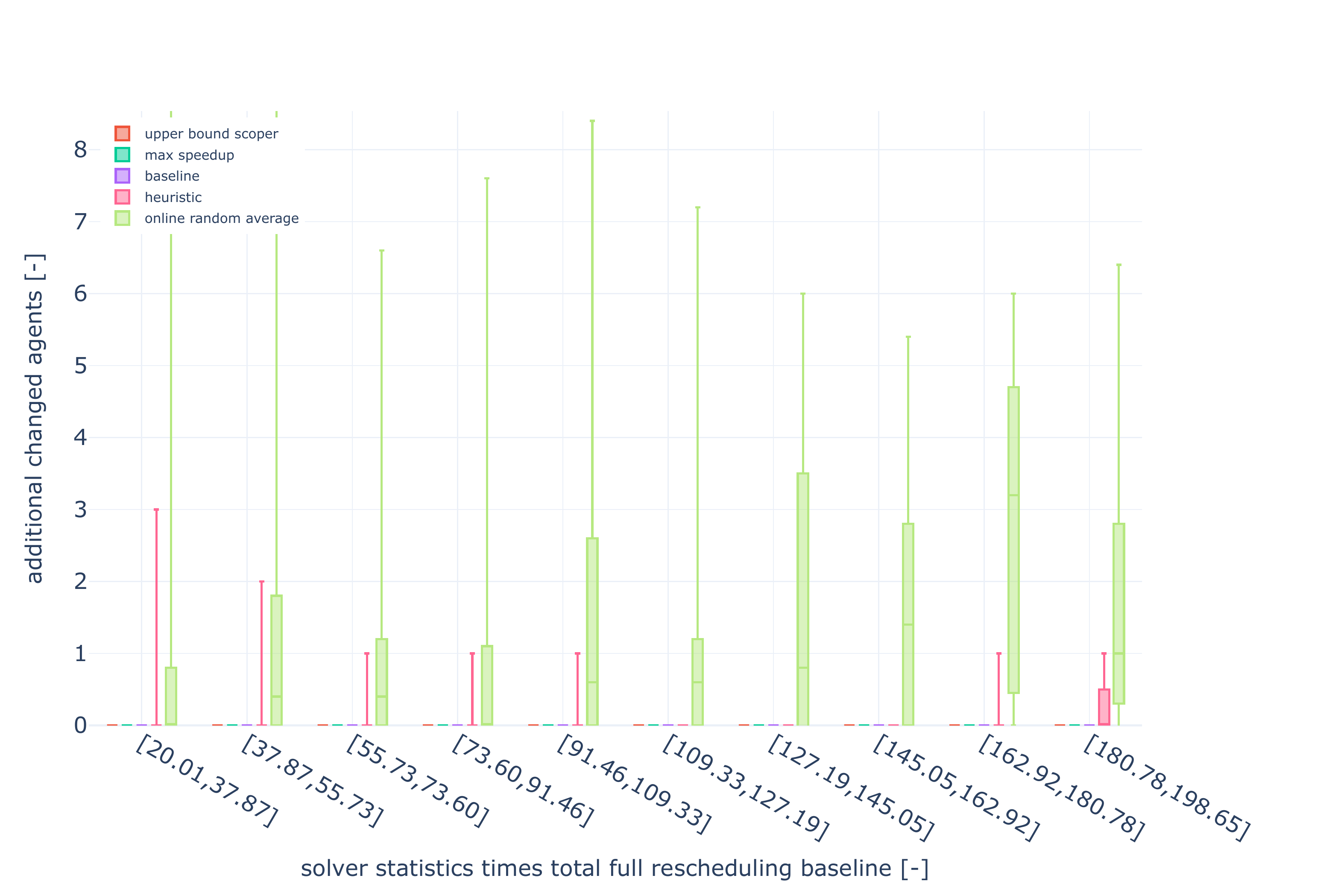}
	\caption{Prediction quality: additional changed trains.}
		\label{fig:additional_changed_trains}
\end{figure}

\subsection{Scope Restriction Quality}
In this section, we present the results about scope restriction quality, i.e. the ability of the scopers to identify the relevant parameters of the re-scheduling problem for a given problem instance. The quality of the scope restriction has a direct impact on the solution time speed-up and the solution quality, which will be discussed in the following sections.

To asses the quality of a scope, we start by defining the \textbf{core problem} as the differences between the original schedule and an optimal re-scheduling solution found with no scope restriction. Theoretically, there could exist multiple core problems with the same optimization score, but here we only focus on the core problem defined by the first optimal solution found by the optimizer.

The quality of any scoper is then defined by its accuracy expressed through the $F_1$-measure. It takes into account the number of trainruns that differ from the core problem, i.e number of false positives (trains not contained in core problem) and false negatives (trains contained in the core problem but not identified by the scoper).

\begin{figure}[hbtp]
	\includegraphics[width=\textwidth]{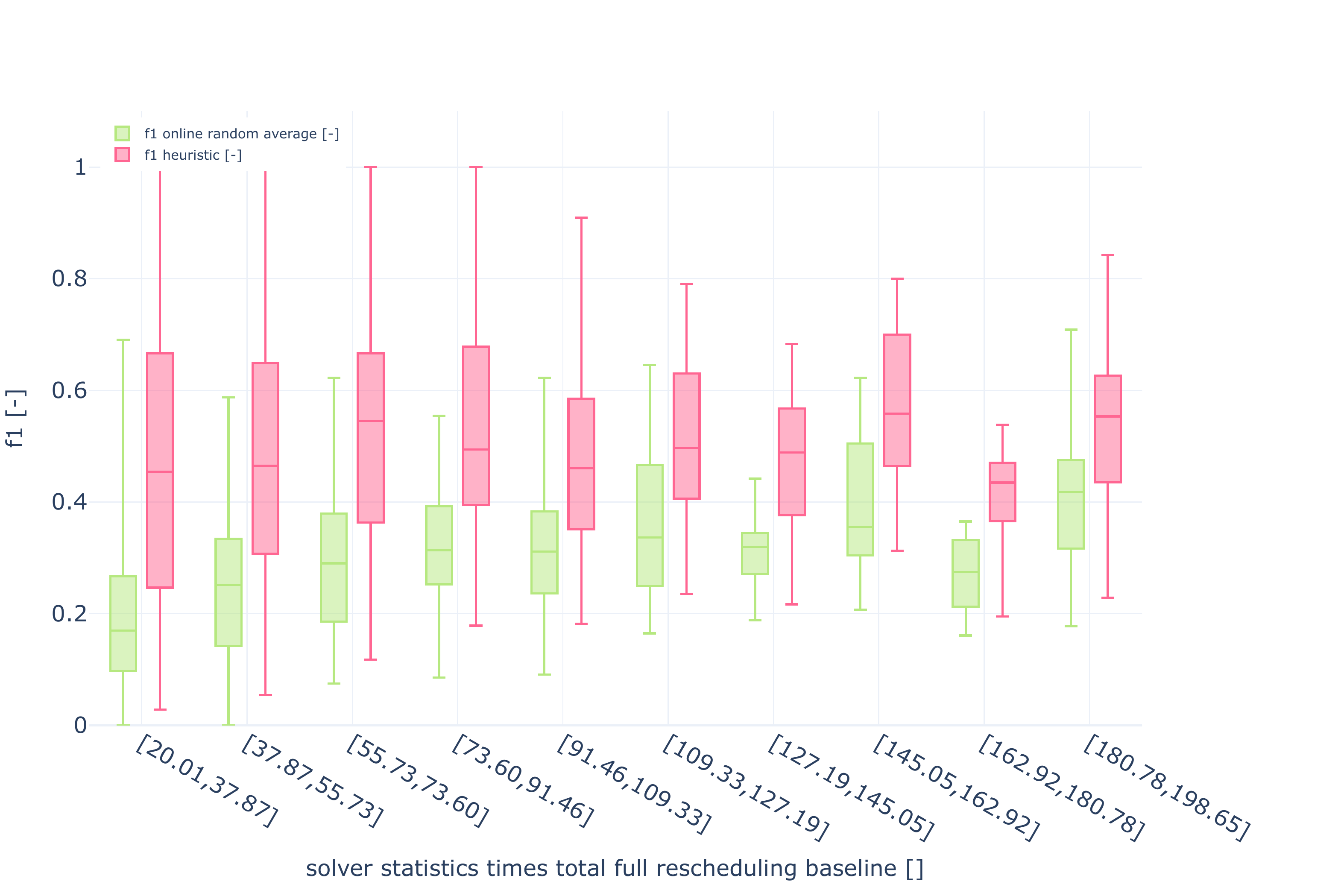}
	\caption{Prediction quality of the heuristic and random scoper. This figure shows the F1 score of the heuristic and random scoper. The performance of the heuristic scoper is well above the random baseline. We lack further baselines to compare with and the heuristic baseline should serve as a benchmark for future improved heuristics or ML scopers.}
	\label{fig:prediction_quality_f1}
\end{figure}

All but the \textbf{heuristic} and \textbf{random} scoper achieve the perfect quality score of $1$, due to their construction. In Figure \ref{fig:prediction_quality_f1} we compare the $F_1$ score of these two scopers to asses their accuracy.

Detailed analysis of the $F_1$ score of the heuristic scoper reveals that there is a large number of false positives and almost no false negatives.

\begin{figure}[hbtp]
	\includegraphics[width=\textwidth]{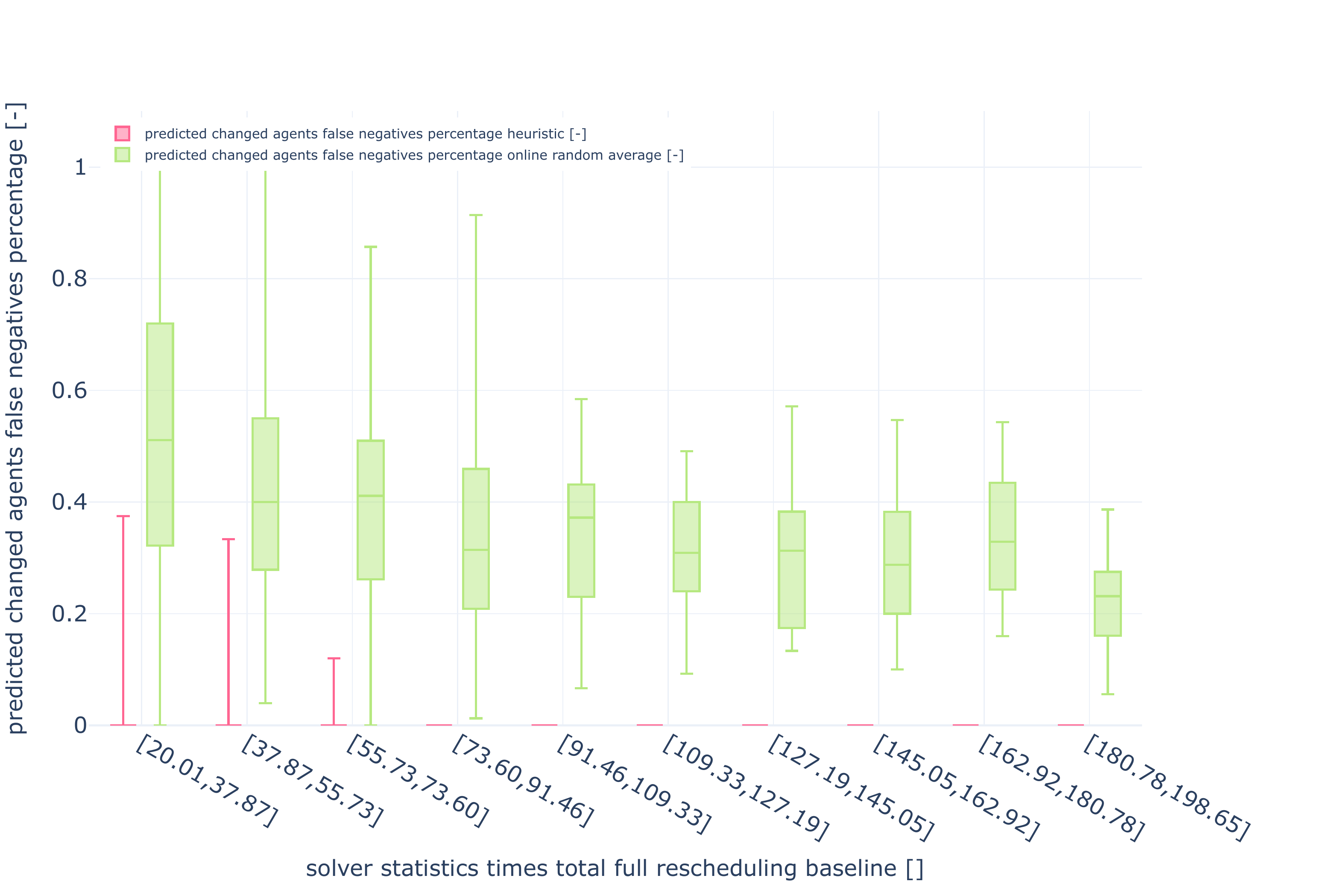}
	\caption{This figure shows the distribution of the prediction quality across all experiments. False negative lead to worse solution quality or infeasible solutions. Only in a few cases did the heuristic scoper produce false negatives, and thus we expect solution quality to be on par with the optimal solution.}
	\label{fig:prediction_quality_false_negative}
\end{figure}

It becomes evident from the number false positives that the heuristic scoper generally overestimates the size of the core problem. The low number of false negatives in the heuristic scope restriction means that the core problem is contained within the predicted scope. This can be extracted from Figure~\ref{fig:prediction_quality_false_negative}.

Given the overestimation of the scope we only expect a minor speed-up for the heuristic scoper while the solution quality should be equivalent to best solution.
In Figure~\ref{fig:prediction_quality_false_positive} and Figure~\ref{fig:additional_changed_trains} in the appendix you find the quality as well as false negatives across all experiments.

The results show that it is difficult to define a heuristic scoper which reliably identifies the core problem without overestimation. This stems from the fact that false negatives lead to poor or even infeasible solutions. We identify this as a key challenge for our approach that needs to be overcome to benefit from the speed-up advantages.

\end{appendices}